\documentclass{article}

\usepackage{stmaryrd,amssymb,amsmath,amsthm,url}

\usepackage{graphicx,version}
\usepackage[usenames]{color}
\newcommand{\NC}{\ensuremath{{\mathbb C}_0}}
\newcommand{\floor}[1]{\lfloor #1 \rfloor}
\newcommand{\ceiling}[1]{\lceil #1 \rceil}
\newcommand{\wortel}[1]{\sqrt[\streep]{#1}}
\newcommand{\streep}{{-}}

\newcommand{\DE}{\axname{DE}}
\newcommand{\SR}{\axname{SquareRoots}}
\newcommand{\FC}{\axname{FC}}
\newcommand{\SA}{\axname{Signs}}
\newcommand{\Mod}{\axname{Mod}}

\newcommand{\NM}{\ensuremath{\mathbb M}}
\newcommand{\sg}{\operatorname{{\mathbf s}}}

\newcommand{\IL}{\axname{IL}}
\newcommand{\IR}{\axname{IR}}
\newcommand{\ILC}{\axname{ILC}}

\newcommand{\Md}{\axname{Md}}

\newcommand{\Ril}{\axname{RIL}}

\newcommand{\NZ}{\ensuremath{\mathbb Z}}
\newcommand{\Nat}{\ensuremath{\mathbb N}}

\newcommand{\axname}[1]{\ensuremath{\textit{#1}}}

\newcommand{\titel}{Cancellation Meadows: a Generic Basis Theorem and Some Applications}

\newtheorem{theorem}{Theorem}
\newtheorem{lemma}{Lemma}
\newtheorem{proposition}{Proposition}
\newtheorem{corollary}{Corollary}
\newtheorem{definition}{Definition}  
\theoremstyle{definition}

\begin{document}

\title{{\titel}~\thanks{ This is a pre-copy-editing, author-produced PDF of an article accepted for publication in The Computer Journal following peer review. The definitive publisher-authenticated version
  [Jan A. Bergstra,
    Inge Bethke,
    and Alban Ponse.
Cancellation Meadows: A Generic Basis Theorem and Some Applications.
The Computer Journal (2013) 56(1): 3-14, 
first published online March 23, 2012, doi:10.1093/comjnl/bxs028] is available online at: 
\texttt{http://comjnl.oxfordjournals.org/content/56/1/}
\texttt{3.full.pdf?keytype=ref\&ijkey=kNeYsWcTYdkTR1u}.
}}

\author{Jan A.\ Bergstra,
    Inge Bethke,
    and
	Alban Ponse
\\[2mm]\small
	  Section Theory of Computer Science\\
	  \small Informatics Institute\\
	  \small University of Amsterdam\\
	  	{\small Url: \url{www.science.uva.nl/~{janb,inge,alban}}}
}

\maketitle

\begin{abstract}
Let $\mathbb{Q}_0$ denote the rational numbers expanded to a ``meadow'', 
that is, after taking its zero-totalized
form ($0^{-1}=0$) as the preferred interpretation. 
In this paper we consider ``cancellation meadows'',
i.e., meadows without proper zero divisors, such as $\mathbb{Q}_0$ 
and prove a 
generic completeness result. 
We apply this result to cancellation meadows expanded with 
differentiation operators, the sign function, and with floor, ceiling and a
signed variant of the square root, respectively.
We give an equational axiomatization of these
operators and thus obtain a finite basis for
various expanded cancellation meadows. 
\\[2mm]
\\[2mm]
{\normalfont\textit{This paper is devoted to the occasion of John Tucker's 60th birthday.
The authors acknowledge his broad scholarly work on algebraic methods in computing. 
In addition Jan Bergstra expresses his great appreciation for over 35 years of joint 
work with John, often unexpectedly emerging from our continuous stream of discussions 
about the field in general.}}
\end{abstract}

\section{Introduction}
\label{sec:Intro}

This paper contributes to the algebraic specification theory of 
number systems.
Advantages and disadvantages of algebraic specification of abstract 
data types
have been amply discussed in the computer science literature and
we do not wish to add anything to those matters here and 
refer to Wirsing~\cite{W90}, 
the seminal 1977-paper~\cite{ADJ} of Goguen \emph{et al.}, 
the overview in Bj{\o}rner and Henson~\cite{CASL}, 
and the ASF+SDF meta-environment
of Klint \emph{et al.}~\cite{Klint2}.

Our focus will be on a particular loose algebraic specification for fields called \emph{meadows},
using the terminology of Broy and Wirsing~\cite{BW81}
who first wrote about loose specifications|i.e.\ the semantic approach not restricted 
to the isomorphism class of initial algebras. 
The theory of algebraic specifications is based on theories of universal algebras.
Some references to universal algebra are,
e.g., Wechler~\cite{Wech} and Graetzer~\cite{Gr}.

The equational specification of the variety of meadows has been proposed 
by Bergstra and Tucker~\cite{BT07} and it has subsequently been
elaborated with more systematic detail in~\cite{BHT}. 
Starting from the signature of fields one obtains the signature of 
meadows by adding a unary inverse operator. At the basis of meadows, 
now, lies the design decision to turn the inverse 
(or division if one prefers a binary notation for pragmatic reasons)
into a total operator by means of the assumption that $0^{-1}=0$. 
By doing so the investigation of number systems as abstract data types 
can be carried out within the original framework of algebraic 
specifications without taking any precautions for partial functions.

Following~\cite{BT07} we write $\mathbb{Q}_0$ for 
the rational numbers expanded to a meadow, that is after 
taking its zero-totalized form as the preferred interpretation. 
The main result
of~\cite{BT07} consists of obtaining an equational initial algebra specification 
of $\mathbb{Q}_0$. The specification takes the form of a general loose specification, valid 
in all fields equipped with a totalized inverse, to which an equation $L_4$ 
specifically designed for the case of rational numbers is taken in 
addition: the equation $L_4$ is based on
Lagrange's theorem that every natural number can be represented as
the sum of 4 squares and reads
\[\frac{1 + x^2 + y^2 + z^2 + u^2}{1 + x^2 + y^2 + z^2 + u^2} = 1.\]
So $L_4$ expresses that for a large collection of numbers $q$,
it holds that $q\cdot q^{-1}=1$
(in particular, those $q$ which can be written as 1 plus the sum of four squares).
Recently, Yoram Hirshfeld has proven that 
\[\frac{1 + x^2 + y^2 }{1 + x^2 + y^2 } = 1\]
suffices (for a proof see~\cite{BM}).

In~\cite{BPvdZ} meadows without proper zero divisors are termed 
\emph{cancellation meadows}. Recently, we found in~\cite{Ono} that meadows were already 
introduced by Komori~\cite{K75} in a report from 1975, where they go by the name of 
\emph{desirable pseudo-fields.}
In~\cite{BHT} it is shown that meadows are precisely the Von Neumann regular rings 
expanded with an inverse operator $\_^{-1}$ and that the equational theory of cancellation meadows 
(there called zero-totalized fields) has a finite basis.  
In this paper we will extend that result to a generic form. 
This enables its application to extended signatures. In particular we 
will examine the case of differential meadows|i.e.\ meadows equipped with differentiation operators. 
A second extension is obtained by adding
a sign function which provides 
one of several mutually interchangeable ways in which the presence of 
an ordering can be equationally specified. The importance of the latter
extension follows from the fact that 
most uses of rational numbers in computer science theory exploit their 
ordering.

We notice that the 
proof of the generic basis theorem is an elaboration of the proof used for the case of closed terms 
that has been dealt with in~\cite{BT07}. 
The proof of the finite basis theorem in~\cite{BHT} uses the 
existence of maximal ideals. Although shorter and simpler, the 
proof via ideals seems not to generalize in the way our proof below 
does. 

Bethke, Rodenburg, and Sevenster~\cite{BR07} demonstrate that finite 
meadows are products of fields, thus 
strengthening  the result in~\cite{BHT} (for the finite case)
that establishes that each meadow can be embedded in a product of 
fields,
a result which was named the \emph{embedding theorem for meadows}. 
We notice
that the basis theorem for meadows, but not its generic form, 
is an immediate consequence of the embedding theorem.

The paper is structured as follows: in the next section we recall
the axioms for meadows and introduce a representation result. 
Then, in Section~\ref{sec:3} we present our main result, the 
generic basis 
theorem. In Section~\ref{sec:4} we introduce differential
meadows. Then, in Section~\ref{sec:5} we extend
cancellation meadows with the sign function. We discuss
a further 
extension 
with floor and ceiling functions and with a square root
in Section~\ref{sec:6}.
 We end the paper in 
Section~\ref{sec:7} with some conclusions.

This paper is compiled from our earlier work as 
reported  in~\cite{BP0,BP,BB}.

\section{Meadows: preliminaries and representation}
\label{sec:2}
In this section we introduce cancellation meadows in detail 
and we discuss
a representation result 
that will be used in Section~\ref{sec:4}.

In~\cite{BT07} \emph{meadows}
were defined as the members of a variety specified by 12 equations. 
However, in~\cite{BHT}
it was established that the 10 equations in Table~\ref{Md}
imply those used in~\cite{BT07}. Summarizing,
a {meadow} is a commutative ring with unit equipped with
a total unary operation $(\_)^{-1}$
named inverse that satisfies the two equations
\begin{align*}
  (x^{-1})^{-1} &= x, \\
  x\cdot(x \cdot x^{-1}) &= x, \quad(\Ril)
\end{align*}
and in which $0^{-1}=0$. Here \Ril\ abbreviates 
\emph{Restricted Inverse Law}.  We write \Md\ for the set of
axioms in Table~\ref{Md}.

\begin{table}
\centering
\hrule
\begin{align*}
	(x+y)+z &= x + (y + z)\\
	x+y     &= y+x\\
	x+0     &= x\\
	x+(-x)  &= 0\\
	(x \cdot y) \cdot  z &= x \cdot  (y \cdot  z)\\
	x \cdot  y &= y \cdot  x\\
	1\cdot x &= x \\
	x \cdot  (y + z) &= x \cdot  y + x \cdot  z\\
	(x^{-1})^{-1} &= x \\
	x \cdot (x \cdot x^{-1}) &= x
\end{align*}
\hrule
\caption{The set \Md\ of axioms for meadows}
\label{Md}
\end{table}

From the axioms in \Md\ the following identities are derivable:
\begin{align*}
	(0)^{-1}  &= 0,\\
	(-x)^{-1} &= -(x^{-1}),\\
	(x \cdot  y)^{-1} &= x^{-1} \cdot  y^{-1},
\\
0\cdot x  &= 0,\\
	x\cdot -y &= -(x\cdot y),\\
	-(-x)     &= x.
\end{align*}

The term \emph{cancellation meadow} is introduced in~\cite{BPvdZ} 
for a zero-totalized field that satisfies the so-called 
``cancellation axiom"
\[
x \neq 0 ~\&~ x\cdot y = x\cdot z ~\longrightarrow~ y=z.\]
An equivalent version of the cancellation axiom that we shall
further use in this paper is the
\emph{Inverse Law} (\IL), i.e., the conditional axiom
\begin{align*}
  x\neq 0 ~\longrightarrow~ x\cdot x^{-1}=1.
 \quad(\IL)
\end{align*}
So \IL\ states that there are no proper zero divisors. 
(Another equivalent formulation of the cancellation property is 
$x\cdot y=0~\longrightarrow~x=0\text{ or }y=0$.)

We write $\Sigma_m=(0,1,+,\cdot,-,^{-1})$ for the
signature of (cancellation) meadows and we shall often write $1/t$ 
or 
\[\frac 1 t \]
for $t^{-1}$, $tu$ for $t\cdot u$, $t/u$ for $t\cdot 1/u$, 
$t-u$ for $t+(-u)$, 
and freely use numerals and exponentiation with constant
integer exponents. 
We shall further write
\[1_x\text{ for } \frac x x\qquad\text{and}\qquad0_x\text{ for }1-1_x,
\]
so, $0_0=1_1=1$, $0_1=1_0=0$, and for all terms $t$,
\[0_t+1_t=1.\]
With the axioms in Table~\ref{Md} we find by \Ril\ that
\begin{align}
\nonumber
1_t\cdot t&=t,\\
\nonumber 
1_t\cdot 1/t&=1/t,\\
\label{eq:1}
(1_t)^2&=1_t,
\end{align}
and we derive the following useful identities:
\begin{align}
\nonumber 
1_t\cdot 0_t&=0,\\
\nonumber 
&\text{(by }1_t\cdot 0_t=1_t(1-1_t)=1_t-1_t=0)\\
\nonumber 
0_t\cdot t&=0,\\
\nonumber 
&\text{(by }(1-1_t)t=t-t=0),\\
\nonumber 
0_t\cdot 1/t&=0\\
\nonumber 
&\text{(by }(1-1_t)1/t=1/t-1/t=0)\\
\label{eq:2}
(0_t)^2&=0_t.\\
\nonumber 
&\text{(by }(1-1_t)^2=1-2\cdot1_t+(1_t)^2=1-1_t=0_t)
\end{align}

In the remainder of this section we discuss a particular standard
representation for meadow terms. We will use this representation
in Section~\ref{sec:4} in order to prove an expressiveness result.

\begin{definition}
A term $P$ over $\Sigma_m$ is a \emph{Standard Meadow Form
(SMF)} if, for some $n\in\Nat$, $P$ is an \emph{SMF of level $n$}.
SMFs of level $n$ are defined as follows:
\begin{description}
\item[\textit{SMF of level $0:$}] each expression of the form $s/t$ 
with $s$ and $t$ ranging over {polynomials}
(i.e., expressions over $\Sigma_m$ without 
inverse operator),
\item[\textit{SMF of level $n+1:$}] each expression of the form
\[0_t\cdot P+1_t\cdot Q\]
with $t$ ranging over {polynomials} and $P$ and $Q$ over SMFs of 
level $n$.
\end{description}
\end{definition}

Observe that if $P$ is an 
SMF of level $n$, then also of level $n+k$ for all $k\in\Nat$.

\begin{lemma} 
\label{La1}
If $P$ and $Q$ are SMFs, then in \Md, 
$P+Q$, $P\cdot Q$, $-P$, and $1/P$
are provably equal to an SMF having the same variables.
\end{lemma}

\begin{proof}
By natural induction on level height $n$.  
We spell out the proof in which \Ril\ is often used.
The mentioned property of having the same variables follows trivially.
\\[2mm]
\textbf{Case $n=0$.}
Let $s,t,u,v$ be polynomials, and $P=s/t$ and $Q=u/v$.

First observe that 
$0_t\cdot s/t=0_t\cdot 1/t\cdot s=0$.
We derive
\begin{align*}
P+Q&=0_t\cdot (P+Q)+1_t\cdot(P+Q)\\
&=0_t\cdot (s/t+u/v)+1_t\cdot (s/t+u/v)\\
&=0_t\cdot u/v+1_t\cdot(s/t+1_t\cdot u/v)
\end{align*}
so it suffices to show that $R=s/t+1_t\cdot u/v$ is equal to an SMF of
level 1:
\begin{align*}
R
&=0_v\cdot (s/t+1_t\cdot u/v)+1_v\cdot(s/t+1_t\cdot u/v)\\
&=0_v\cdot s/t+1_v\cdot(s/t\cdot 1_v+1_t\cdot u/v)\\
&=0_v\cdot s/t+1_v\cdot(\frac{sv+tu}{tv}).
\end{align*}

The remaining cases are trivial:
\[
P\cdot Q=su/tv,\quad
-P=-s/t,\quad\text{and}\quad
1/P=t/s.
\]
\textbf{Case $n+1$.}
Let $P=0_t\cdot S+1_t\cdot T$ and $Q=0_s\cdot U+1_s\cdot V$ 
with $S,T,U,V$ all SMFs of level $n$. 

We first derive
\begin{align*}
P+Q&=0_t\cdot P+1_t\cdot P+Q\\
&=0_t\cdot (S+Q)+1_t\cdot(T+Q)\\ 
&=0_t\cdot (0_s\cdot(S+U)+1_s\cdot(S+V))~+\\
&\phantom{~=}
1_t\cdot(0_s\cdot(T+U)+1_s\cdot(T+V))
\end{align*}
and by induction each of the pairwise sums of $S,T,U,V$ 
equals some SMF.

Next, we derive
\begin{align*}
P\cdot Q&=0_s\cdot P \cdot U+1_s\cdot P\cdot V\\
&=0_s\cdot(0_t\cdot S\cdot U+1_t\cdot T\cdot U)~+\\
&\phantom{~=}1_s\cdot(0_t\cdot S\cdot V+1_t\cdot T\cdot V)
\end{align*}
and by induction each of
the pairwise products of $S,T,U,V$ 
equals some SMF.

Furthermore, $-P=0_t\cdot(-S) + 1_t\cdot (-T)$, 
which by induction is provably equal to an SMF.

Finally, $1/P=0_t\cdot(1/P)+1_t\cdot(1/P)$, hence
\begin{align*}
1/P&=0_t\cdot \frac 1 {0_t\cdot S+1_t\cdot T}+1_t\cdot 
\frac 1 {0_t\cdot S+1_t\cdot T}\\
&=0_t\cdot \frac{0_t}{0_t\cdot(0_t\cdot S+1_t\cdot T)}+
1_t\cdot \frac{1_t} {1_t\cdot(0_t\cdot S+1_t\cdot T)}\\
&=0_t\cdot \frac{0_t}{0_t\cdot S}+1_t\cdot \frac{1_t} {1_t\cdot T}\\
&=0_t\cdot {1}/{S}+1_t\cdot {1}/{T}
\end{align*}
and by induction there exist SMFs $S'$ and $T'$ such that $S'=1/S$ and 
$T'=1/T$, hence $1/P=0_t\cdot S'+1_t\cdot T'$.
\end{proof}

\begin{theorem}
\label{SMF}
Each term over $\Sigma_m$ can be represented by an SMF 
with the same variables.
\end{theorem}

\begin{proof}
By structural induction. Let $P$ be a term over $\Sigma_m$.
If $P= 0$ or $P=1$ or $P=x$, then $P=P/1$, and the latter is an SMF of level 0.
 The other cases follow immediately from Lemma~\ref{La1}.
\end{proof}

\section{A generic  basis theorem}
\label{sec:3}
In this section we prove a finite basis result
for the equational theory of cancellation meadows. 
This result is formulated in a generic way so that
it can be used for any expansion of a meadow that satisfies
the propagation properties defined below.

\begin{definition}
Let $\Sigma$ be an extension of $\Sigma_m=(0,1,+,\cdot,-,^{-1})$, the
signature of meadows. Let $E\supseteq \Md$ (with \Md\ the set
of axioms for meadows given in Table~\ref{Md}).
\begin{enumerate}
\item
$(\Sigma,E)$ has the \emph{propagation property for pseudo units} if for
each pair of $\Sigma$-terms $t,r$ and context $C[~]$,
\[E\vdash 1_t\cdot C[r]=1_t\cdot C[1_t\cdot r].\]
\item
$(\Sigma,E)$ has the \emph{propagation property for pseudo zeros} if for
each pair of $\Sigma$-terms $t,r$ and context $C[~]$,
\[E\vdash 0_t\cdot C[r]=0_t\cdot C[0_t\cdot r].\]
\end{enumerate}
\end{definition}
Preservation of these propagation properties admits the following 
nice result:
\begin{theorem}
[Generic Basis Theorem for Cancellation Meadows]
\label{st}
If $\Sigma\supseteq \Sigma_m,~ E\supseteq \Md$ and $(\Sigma,E)$ 
has the
pseudo unit propagation property and the pseudo zero propagation property,
then $E$ is a basis (a complete axiomatization) of 
$\Mod_\Sigma(E\cup\IL)$.
\end{theorem}
The structure of our proof of this theorem is as follows: 
let $r=r(\overline x)$ and $s=s(\overline x)$ be $\Sigma$-terms and 
let 
$\overline c$ be a series of fresh constants. We write
$\Sigma(\overline c)$ 
for the signature extended with these constants.  Then
\begin{align}
\nonumber
 &E\cup\IL \models r=s \quad\text{in }\Sigma\\
\label{a} 
 &\iff 
 E\cup\ILC\models r(\overline c)=s(\overline c)
\quad\text{in }\Sigma(\overline c)
\\
\label{e}
&\iff  E \vdash_{\IR} r(\overline c)=s(\overline c)\\
\label{f}
&\iff
 E \vdash r=s\quad\text{in }\Sigma.
\end{align}
Here provability ($\vdash$) refers to equational logic;
the notation further used means this:
\begin{itemize}
\item \ILC, the \emph{Inverse Law for Closed terms} is the set
$\{t=0\vee 1_t=1\mid t\in T(\Sigma(\overline c))\}$,
where $T(\Sigma(\overline c))$ denotes the set of 
closed terms over $\Sigma(\overline c)$.
\item \IR\ is the \emph{Inverse Rule}:
$E\vdash_{\IR} r=s$ means that $\exists k\in \Nat$ s.t.
$E\vdash_{\IR}^k r=s$, and
$E\vdash_{\IR}^k r=s$ means that
$E\vdash r=s$ provided that the rule
\[
\IR\qquad
\frac{E\cup\{t=0\}\vdash r=s\quad E\cup\{1_t=1\}\vdash r=s}{E\vdash r=s}
\]
with $t$ ranging over
$T(\Sigma(\overline c))$ may be used $k$ times.
\end{itemize}
Before we prove Theorem~\ref{st} | i.e., equivalences 
\eqref{a}--\eqref{f} | we establish the following
preliminary result:

\begin{proposition} 
\label{A}
Assume $\Sigma\supseteq \Sigma_m$,
$E\supseteq\Md$ and $(\Sigma,E)$ has the
propagation property for pseudo units and for pseudo zeros. 
Then for $t,r,s\in T(\Sigma)$,
\begin{eqnarray}
\label{g}
E\cup\{t=0\}\vdash_{\IR} r=s &\Longrightarrow&
E\vdash 0_t\cdot r=0_t\cdot s,\\
\label{h}
E\cup\{1_t=1\}\vdash_{\IR} r=s &\Longrightarrow&
E\vdash 1_t\cdot r=1_t\cdot s.
 \end{eqnarray}
\end{proposition}

\begin{proof}
We prove 
\begin{eqnarray}
\label{i}
E\cup\{t=0\}\vdash_{\IR}^k r=s &\Longrightarrow&
E\vdash 0_t\cdot r=0_t\cdot s,\\
\label{j}
E\cup\{1_t=1\}\vdash_{\IR}^k r=s &\Longrightarrow&
E\vdash 1_t\cdot r=1_t\cdot s
\end{eqnarray}
simultaneously by induction on $k$. We use the symbol $\equiv$
to denote syntactic equivalence.
\begin{description}
\item[Case $k=0$.] By induction on proof lengths. 
For \eqref{i} the  only interesting case is $(r=s)\equiv (t=0)$, 
so we have to show that $E\vdash 0_t\cdot t=0_t\cdot 0$. This
follows directly from $E\supseteq\Md$.

For \eqref{j} the only interesting case is $(r=s)\equiv (1_t=1)$,
and also $E\vdash (1_t)^2=1_t\cdot 1$
follows directly from $E\supseteq\Md$.
\item[Case $k+1$.] By induction on the length 
of the proofs of 
$E\cup\{t=0\}\vdash_{\IR}^{k+1} r=s$ and 
$E\cup\{1_t=1\}\vdash_{\IR}^{k+1} r=s$. 
There are 3 interesting cases for each of \eqref{i} and \eqref{j}:
\begin{enumerate}
\item The $\vdash_{\IR}^{k+1}$ results follow from the assumption
$(r=s)\equiv (t=0)$ or $(r=s)\equiv (1_t=1)$, respectively. 
These results follow in the same way as above.
\item The $\vdash_{\IR}^{k+1}$ results follow from the context rule, 
so $r\equiv C[v]$, $s\equiv C[w]$ and
\begin{enumerate}
\item[\eqref{i}] $E\cup\{t=0\}\vdash_{\IR}^{k+1} v=w$.
By induction, $E\vdash 0_t\cdot v=0_t\cdot w$. Hence, 
$E\vdash 0_t\cdot C[0_t\cdot v]=0_t\cdot C[0_t\cdot w]$, 
and by $(\Sigma,E)$ having the propagation property for
pseudo zeros, $E\vdash 0_t\cdot C[v]=0_t\cdot C[w]$.
\item[\eqref{j}] $E\cup\{1_t=1\}\vdash_{\IR}^{k+1} v=w$.
By induction, $E\vdash 1_t\cdot v=1_t\cdot w$. Hence, 
$E\vdash 1_t\cdot C[1_t\cdot v]=1_t\cdot C[1_t\cdot w]$, 
and by $(\Sigma,E)$ having the propagation property for
pseudo units, $E\vdash 1_t\cdot C[v]=1_t\cdot C[w]$.
\end{enumerate}
\item The $\vdash_{\IR}^{k+1}$ results follow from the \IR\ rule, that is
\begin{enumerate}
\item[\eqref{i}] $E\cup\{t=0\}\cup\{h=0\}\vdash_{\IR}^{k} r=s$
and $E\cup\{t=0\}\cup\{1_h=1\}\vdash_{\IR}^{k} r=s$.
By induction, $E \cup\{h=0\}\vdash 0_t\cdot r=0_t\cdot s$ and 
$E \cup\{1_h=1\}\vdash 0_t\cdot r=0_t\cdot s$. 
Again applying induction ($\vdash$ derivability implies $\vdash_{\IR}^k$ 
derivability) 
yields
\begin{align*}
&E\vdash 0_h\cdot 0_t\cdot r=0_h\cdot 0_t\cdot s,\\
&E\vdash 1_h\cdot 0_t\cdot r=1_h\cdot 0_t\cdot s.
\end{align*} 
We derive $0_t\cdot r=(0_h+1_h)\cdot0_t\cdot r=
0_h\cdot 0_t\cdot r+1_h\cdot0_t\cdot r=0_h\cdot 0_t\cdot s+1_h\cdot 
0_t\cdot s=0_t\cdot s$.
\item[\eqref{j}] $E\cup\{1_t=1\}\cup\{h=0\}\vdash_{\IR}^{k} r=s$
and $E\cup\{1_t=1\}\cup\{1_h=1\}\vdash_{\IR}^{k} r=s$.
Similar.
\end{enumerate}
\end{enumerate}
\end{description}
\end{proof}

\begin{proof}[Proof of Theorem~\ref{st}.]
We now give a detailed proof of equivalences \eqref{a}--\eqref{f}, using Proposition~\ref{A}. 
For model theoretic details we refer 
to~\cite{CK90}.
\begin{itemize}

\item[\eqref{a}]
($\Longrightarrow$)
Assume $E\cup \IL\models r=s$.
Let $\NM$ be a model of $E\cup \ILC$ 
(over $\Sigma(\overline c)$). 
Then $\NM\models r(\overline c)=s(\overline c)$ if and only 
if $\NM'\models r(\overline c)=s(\overline c)$ for $\NM'$ the minimal submodel of $\NM$. Now
$\NM'$ is also a model for \IL\ because \ILC\ concerns all closed terms 
and each value in the domain of $\NM'$ is the interpretation
of a closed term. 
So, by assumption $\NM'\models r=s$,
and, in particular (by substitution),
$\NM'\models r(\overline c)=s(\overline c)$.

($\Longleftarrow$) Assume $E\cup \ILC\models 
r(\overline  c)=s(\overline  c)$.
Let $\NM$ be a model of $E\cup \IL$ (over $\Sigma$). We have to 
show that $\NM\models r(\overline x)=s(\overline x)$, or, stated differently, that 
for $\overline  a=a_1,...,a_n$ a series of
values from $\NM$'s domain,
$(\NM,x_i\mapsto a_i)\models r = s$ where $x_i\mapsto a_i$ represents the 
assignment of $a_i$ to $x_i$. 
Extend $\Sigma$ with a fresh constant $c_i$ for
each $a_i$ and let $\NM(\overline  c)$ be the expansion 
of $\NM$ in which each constant $c_i$ is interpreted as $a_i$. 
Then $\NM(\overline  c)$
satisfies \ILC\ because $\NM$ satisfies \IL, so by assumption
$\NM(\overline  c)\models r(\overline c)=s(\overline c)$, and therefore
$(\NM(\overline  c),x_i\mapsto a_i)\models r=s$ and thus also
$(\NM,x_i\mapsto a_i)\models r=s$, as was to be shown.

\item[\eqref{e}]
($\Longrightarrow$) Let $E^C$ be the set of all closed
instances over the extended signature $\Sigma(\overline c)$, then 
\[E^C\cup\ILC\models r(\overline c)=s(\overline c).\]
By compactness there is a finite set 
$F\subseteq E^C\cup\ILC$ such that 
$F\models r(\overline c)=s(\overline c).$

Now apply induction on the number of elements from \ILC\ in $F$, say $k$.
\begin{description}
\item[Case $k=0$.] By completeness we find 
$E \vdash r(\overline c)=s(\overline c)$, and thus 
$E \vdash_{\IR} r(\overline c)=s(\overline c)$.
\item[Case $k+1$.] Assume $(t=0\vee 1_t=1)\in F$ and let 
$F'=F\setminus \{t=0\vee 1_t=1\}$. 
Reasoning in propositional logic we find 
\[F' \models  (t=0\vee 1_t=1)\rightarrow
r(\overline c)=s(\overline c)\]
and thus 
\begin{align*}
F' \models  &(t=0\rightarrow r(\overline c)=s(\overline c))
~\wedge\\&(1_t=1\rightarrow
r(\overline c)=s(\overline c)),
\end{align*}
which in turn is
equivalent with 
\begin{align*}
F'\cup \{t=0\}\models 
r(\overline c)=s(\overline c),\\
F'\cup \{1_t=1\}\models 
r(\overline c)=s(\overline c).
\end{align*}
By induction, 
$E\cup\{t=0\}\vdash_{\IR} r(\overline c)=s(\overline c)$ and
$E\cup\{1_t=1\}\vdash_{\IR} r(\overline c)=s(\overline c)$, 
and thus by \IR,
\[E \vdash_{\IR} r(\overline c)=s(\overline c).\]
\end{description}
($\Longleftarrow$) This follows from the
soundness of \IR\ with respect to \ILC. That is, if
$E\vdash u=v$ because
$E\cup\{t=0\}\vdash u=v$ and $E\cup\{1_t=1\}\vdash u=v$, then $E
\cup\{t=0\vee 1_t=1\} \models u=v$, so $E\cup\ILC\models u=v$.

\item[\eqref{f}]

($\Longrightarrow$) By induction on the length of the proof, 
using Proposition~\ref{A}: if $E\vdash_{\IR} r(\overline c)=s(\overline c)$ 
follows from \IR\ (the only interesting case), then 
\begin{align*}
E\cup\{t=0\}\vdash_{\IR} r(\overline c)=s(\overline c),\\
E\cup\{1_t=1\}\vdash_{\IR} r(\overline c)=s(\overline c),
\end{align*}
so  $E\vdash 0_t\cdot r(\overline c)=0_t\cdot s(\overline c)$ 
by~\eqref{g} and 
$E\vdash 1_t\cdot r(\overline c)=1_t\cdot s(\overline c)$ 
by~\eqref{h}. Thus 
\begin{align*}
E\vdash r(\overline c)=(0_t+1_t)\cdot r(\overline c)&=0_t\cdot 
r(\overline c)+1_t\cdot r(\overline c)\\
&=0_t\cdot s(\overline c)+1_t\cdot s(\overline c)\\
&=s(\overline c).
\end{align*}
A similar proof result is obtained by replacing $r(\overline c)$ 
by $r$ and $s(\overline c)$ by $s$.

($\Longleftarrow$) Trivial: if $E\vdash r=s$, then 
$E\vdash r(\overline c)=s(\overline c)$ in the extended 
signature $\Sigma(\overline c)$.
So, $E\vdash_{\IR} r(\overline c)=s(\overline c)$.

\end{itemize}
\end{proof}

A first application of Theorem~\ref{st} concerns the equational theory of cancellation meadows:

\begin{corollary}
The set of axioms $\Md$ (see Table~\ref{Md})
is a finite basis (a complete axiomatization) of 
$\Mod_{\Sigma_m}(\Md\cup\IL)$. 
\end{corollary}

\begin{proof}
It remains to be
shown that the propagation properties for pseudo units 
and for pseudo zeros hold in \Md. This follows easily
by case distinction on the forms that $C[r]$
may take and the various identities on $1_t$ and $0_t$.
As an example consider the case
$C[\_]\equiv \_+u$. Then 
\begin{align*}
1_t\cdot C[r]&=1_t\cdot (r+u)\\
&=1_t\cdot r + 1_t\cdot u\\
&=1_t\cdot 1_t\cdot r+1_t\cdot u\\
&=1_t\cdot C[1_t\cdot r].
\end{align*}
The remaining cases can be proved in a similar way.
\end{proof}

\section{Differential Meadows}
\label{sec:4}
In this section
we provide an elegant equational axiomatization of differential operators
and with the generic basis theorem we obtain a finite basis for differential cancellation meadows.

\subsection{Differential Meadows}
Given some $n\geq 1$ we extend the signature $\Sigma_m$ of meadows 
with differentiation operators and constants $X_1,...,X_n$ to model
functions to be differentiated:
\[\frac{\partial}{\partial X_i}: \NM\rightarrow \NM\]
for $i=1,...,n$ and some meadow \NM. 
We write $\Sigma_{md}$ for this extended signature.
Equational axioms for $\frac{\partial}{\partial X_i}$
are given in Table~\ref{t1}, where \eqref{D4} and 
\eqref{D5} define $n^2$ 
equational axioms. 
Observe that the \Md\ axioms together with Axiom~\eqref{D3} imply
\[
\frac{\partial}{\partial X_i}(0)=0. 
\]
Furthermore, using Axiom~\eqref{D1} one easily proves: 
\[
\frac{\partial}{\partial X_i}(-x) =-\frac{\partial}{\partial X_i}(x).
\]

\begin{table}
\hrule
\begin{align}
\label{D1}
\frac{\partial}{\partial X_i}(x+y)&=\frac{\partial}{\partial X_i}(x)
+\frac{\partial}{\partial X_i}(y)\\
\label{D2}
\frac{\partial}{\partial X_i}(x\cdot y)&=
\frac{\partial}{\partial X_i}
(x)\cdot y+x\cdot \frac{\partial}{\partial X_i}(y)\\
\label{D3}
\frac{\partial}{\partial X_i}(x\cdot x^{-1})&=0\\
\label{D4}
\frac{\partial}{\partial X_i}(X_i)&=1\\
\label{D5}
\frac{\partial}{\partial X_i}(X_j)&=0 \quad
\text{ if } i \neq j
\end{align}
\hrule
\caption{The set \DE\ of axioms for differentiation}
\label{t1}
\end{table}

First we establish the expected corollary of Theorem~\ref{st}:
\begin{corollary}
\label{cor:diff}
The set of axioms $\Md\cup\DE$ (see Tables~\ref{Md} and \ref{t1})
is a complete axiomatization of 
$\Mod_{\Sigma_{md}}(\Md\cup\DE\cup\IL)$.
\end{corollary}
\begin{proof}
The pseudo unit propagation property requires a check for 
$\frac{\partial}{\partial X_i}(\_)$ only:
\begin{equation}
\label{tja}
\frac{\partial}{\partial X_i}(1_t\cdot r)=
\frac{\partial}{\partial X_i}(1_t)\cdot r+1_t\cdot 
\frac{\partial}{\partial X_i}(r)=
1_t\cdot\frac{\partial}{\partial X_i}(r).
\end{equation}
Multiplication with $1_t$ now yields the property.
From \eqref{tja} we get
\begin{align*}
0_t\cdot \frac{\partial}{\partial X_i}(r)
&=\frac{\partial}{\partial X_i}(r)-1_t\cdot 
\frac{\partial}{\partial X_i}(r)\\
&\stackrel{\eqref{tja}}=\frac{\partial}{\partial X_i}(r)-
\frac{\partial}{\partial X_i}(1_t\cdot r)
=
\frac{\partial}{\partial X_i}(0_t\cdot r)
\end{align*}
and multiplication with $0_t$
then yields the pseudo zero propagation property.
\end{proof}

A \emph{differential meadow} is a meadow
equipped with formal variables $X_1,...,X_n$ and differentiation
operators $\frac{\partial}{\partial X_i}(\_)$ that satisfies 
the axioms in \DE.

We conclude this section with an elegant consequence of the 
fact that we are working in the setting of meadows, namely 
the consequence that the differential of an inverse follows 
from the \DE\ axioms.

\begin{proposition}
\[
\Md\cup\DE\,\vdash \frac{\partial}{\partial X_i}(1/x)=-(1/x^2) 
\cdot 
\frac{\partial}{\partial X_i}(x).\]
\end{proposition}

\begin{proof}
By Axioms~\eqref{D3} and~\eqref{D2},
\[
0=\frac{\partial}{\partial X_i}(x/x)
=
\frac{\partial}{\partial X_i} (x)\cdot 1/x
+x\cdot\frac{\partial}{\partial X_i}(1/x),
\]
so 
\begin{align*}
0&=0\cdot(1/x)= \displaystyle
\frac{\partial}{\partial X_i}(x/x)\cdot(1/x)\\
&=\frac{\partial}{\partial X_i}(x)\cdot 1/x^2+(x/x)\cdot
\frac{\partial}{\partial X_i}(1/x)\\
& \displaystyle\stackrel{\eqref{tja}}=1/x^2 \cdot 
\frac{\partial}{\partial X_i}(x)
+
\frac{\partial}{\partial X_i}((x/x)\cdot(1/x))\\
&\stackrel{\Ril}=1/x^2 \cdot \frac{\partial}{\partial X_i}(x)+
\frac{\partial}{\partial X_i}(1/x),
\end{align*}
and hence
\[
\frac{\partial}{\partial X_i}(1/x)=-(1/x^2) \cdot 
\frac{\partial}{\partial X_i}(x).
\]
\end{proof}

\subsection{Existence of Differential Meadows}
\label{subsec:4.2}
In this section we show the \emph{existence} of differential 
meadows with formal variables $X_1,...,X_n$ for arbitrary 
finite $n>0$. 
First we define a particular cancellation meadow, and then we expand 
this meadow to a differential cancellation meadow by adding formal
differentiation.

\paragraph{The Zariski topology congruence over $\NC^n$.}
We will use some terminology from algebraic geometry, 
in particular we will use the Zariski topology~\cite{Zar44,Har77}. 
Open (closed) sets in this topology will be indicated 
as Z-open (Z-closed). Recall that 
complements of Z-closed sets are Z-open and complements of 
Z-open sets are Z-closed,
finite unions of Z-closed sets are Z-closed, and
intersections of Z-closed sets are Z-closed.
Let $\NC$ denote the zero-totalized expansion of the complex numbers. 
We will make use of the following facts:

\begin{enumerate} 
\item \label{een}
The solutions of a set of polynomial equations (with $n$ or less 
variables) within $\NC^n$ constitute a Z-closed subset of $\NC^n$. 
Here 'polynomial' has the conventional meaning, not involving division. 
Taking equations $1=0$ and $0=0$ respectively, it follows that both 
$\emptyset$ and $\NC^n$ are Z-closed (and Z-open as well).

\item \label{twee}
Intersections of non-empty Z-open sets are non-empty.
\end{enumerate}
In the following we consider terms 
\[t(\overline X)=t(X_1,...,X_n)\]
with $t=t(\overline x)$ a $\Sigma_m$-term and we write 
$T(\Sigma_m(\overline X))$ for the set of these terms.
For $V\subseteq\NC^n$ we define the equivalence
\[\equiv^V_{\NC^n}\]
on $T(\Sigma_m(\overline X))$
by $t(\overline X)\equiv^V_{\NC^n}r(\overline X)$ 
if each assignment $\overline X\mapsto V$
evaluates both sides to  equal values in \NC. 
It follows immediately that for each $V\subseteq\NC^n$, 
$T(\Sigma_m(\overline X))/\equiv^V_{\NC^n}$ is a 
meadow. In particular, if 
$ V = \emptyset$ one obtains the trivial meadow ($0=1$) as both 0 
and 1 satisfy any universal quantification over an empty set. 
If $V$ is a singleton this quotient is a cancellation meadow. 
In other cases the meadow may not satisfy 
the cancellation property. 
Indeed, suppose that $n=1$ and $V = \{0,1\}$ and
let $t(X) = X$. Now $t(1) \neq 0$. Thus
$t(X) \neq 0$ in  $T(\Sigma_m(X))/\equiv^V_{\NC}$. 
If that is assumed to be a cancellation meadow, however, 
one has $1_{t(X)}=1$, but $1_{t(0)}=0$, thus refuting 
$1_{t(X)}=1$.

We now define the relation $\equiv_{ZTC}$ 
(Zariski Topology Congruence over $\NC^n$) by
\[t\equiv_{ZTC} r \iff \exists V (
  \text{$V$ is Z-open, $V\neq\emptyset$ and 
        $t\equiv^V_{\NC^n}r$)}.
\]
The relation $\equiv_{ZTC}$ is indeed a congruence 
for all meadow operators: 
the equivalence properties follow easily; 
for $0\equiv_{ZTC} 0$ and $1\equiv_{ZTC} 1$,
take $V=\NC^n$, and if 
$P\equiv_{ZTC} P'$ and $Q\equiv_{ZTC} Q'$, witnessed 
respectively by $V$ and $V'$, 
then 
\[P+P'\equiv_{ZTC} Q+Q'\quad\text{and}\quad
P\cdot P'\equiv_{ZTC} Q\cdot Q'\] 
are witnessed by
$V\cap V'$  which is Z-open and non-empty because  of
fact~\ref{twee} above. Finally $-P\equiv_{ZTC} -P'$ and 
$(P)^{-1}\equiv_{ZTC}(P')^{-1}$
are both witnessed by $V$.

Theorem~\ref{SMF}, i.e., the (SMF) representation result for meadow terms implies for 
\[T(\Sigma_m(\overline X))/\equiv_{ZTC}\]
that each term can be 
represented by 0 or by $p/q$ with $p$ and $q$ polynomials not 
equal to 0. We notice that it is decidable whether or not a 
polynomial equals the 0-polynomial by taking all corresponding 
products of powers of the $X_1,...,X_n$ together and then 
checking that all coefficients vanish. 

As an example, let $P$ be the SMF of level 1 defined by
\[P=0_{1-X_1}\cdot \frac{2X_1}{X_2}+1_{1-X_1}\cdot 
\frac{1+X_2-2X_1X_3}{8-X_1X_3^2}.\]
Now in $T(\Sigma_m(\overline X))/\equiv_{ZTC}$, the polynomial 
$1-X_1$ is on some Z-open non-empty set $V$ not equal to 0 
(see fact~\ref{een} above), thus
$1_{1-X_1}\equiv^V_{\NC^n}1$ and $0_{1-X_1} \equiv^V_{\NC^n} 0$, 
and hence 
\[P\equiv_{ZTC}\frac{1+X_2-2X_1X_3}{8-X_1X_3^2}.\]
So, in $T(\Sigma_m(\overline X))/\equiv_{ZTC}$, 
the SMF level-hierarchy collapses and
terms can be represented by either $0$ or by $p/q$ with both $p$ and 
$q$ polynomials not equal to 0. In the second case
$1_{p/q} = 1$ and therefore it is a cancellation meadow.
Furthermore, equality is decidable in this model. 
Indeed to check that $1_p = 1$ (and $0_p = 0$) for a polynomial 
$p$ it suffices to check that $p$ is not 0 over the complex numbers. 
Using the SMF representation all closed terms are either 0 or 
take the form $p/q$ with $p$ and $q$ nonzero polynomials. 
For $q$ and $q'$ nonzero 
polynomials we find that 
\[p/q \equiv_{ZTC} p'/q' 
\iff p\cdot q' - p'\cdot q = 0\]
which we have already found to be decidable.

\paragraph{Constructing a differential cancellation meadow.}
In $T(\Sigma_m(\overline X))/\equiv_{ZTC}$ the differential 
operators can be defined as follows:
\[\frac{\partial}{\partial X_i}(0)=0\]
and, using the fact that differentials on polynomials are known,
\[\frac{\partial}{\partial X_i}(\frac p q) =
\frac{\frac{\partial}{\partial X_i}(p) \cdot 
q - p \cdot \frac{\partial}{\partial X_i}(q)}{q^2}.\]
Let $V$ be the set of 0-points of $q$ and let $U={\sim} V$, 
the complement of $V$. 
Then $p/q$ is  differentiable on $U$ and the derivative coincides 
with the formal derivative used in the definition. 
This definition is representation independent:
consider $p'/q'\equiv_{ZTC}p/q$ with $V'$ the 0-points of $q'$
and $U'={\sim}V'$. Then there is some 
non-empty and Z-open $W$
such that $p/q\equiv^W_{\NC^n}p'/q'$. Now
$W\cap U\cap U'$ is non-empty and Z-open, and on this set,
\[\frac{\partial}{\partial X_i}(\frac p q)=
\frac{\partial}{\partial X_i}(\frac {p'}{q'}).\]
So, formal differentiation $\partial/\partial X_i$ preserves 
the congruence properties. Finally, we check the soundness of the
\DE\ axioms:
\\[2mm]
Axiom~\eqref{D1}:
Consider $t+t'$. In the case that one of $t$ 
and $t'$ equals 0, axiom D1 is obviously sound. 
In the remaining case,
$t=p/q$ and $t'=p'/q'$ with all polynomials
not equal to 0 and 
\[t+t'=\frac{pq'+p'q}{qq'}. 
\]
Using ordinary differentiation on polynomials we derive
\begin{align*}
&{\frac{\partial}{\partial X_i}(t+t')}\\
&= \frac{\frac{\partial}{\partial X_i}(pq'+p'q) \cdot 
qq' - (pq'+p'q) \cdot \frac{\partial}{\partial X_i}(qq')}{(qq')^2}\\
&=\frac{\frac{\partial}{\partial X_i}(p)\cdot q\cdot (q')^2+
\frac{\partial}{\partial X_i}(p')\cdot q^2\cdot q' }{(qq')^2}~+\\
&\phantom{~=}
\frac{-p\cdot \frac{\partial}{\partial X_i}(q)\cdot (q')^2
- p'\cdot \frac{\partial}{\partial X_i}(q')\cdot q^2}{(qq')^2}\\
&=\frac{\partial}{\partial X_i}(\frac{p}{q})\cdot 1_{(q')^2}+
\frac{\partial}{\partial X_i}(\frac{p'}{q'})\cdot 1_{q^2}\\
&=\frac{\partial}{\partial X_i}(t)+\frac{\partial}{\partial X_i}(t').
\end{align*}
\\[2mm]
Axiom~\eqref{D2}: Similar.
\\[2mm]
Axiom~\eqref{D3}:
Consider $t$, then either $t=0$ or $t/t=1$, and in both cases
$\displaystyle\frac{\partial}{\partial X_i}(\frac t t)=0$.
\\[2mm]
Axioms schemes~\eqref{D4} and~\eqref{D5}: We derive
\[\frac{\partial}{\partial X_i}(X_j)=
\frac{\partial}{\partial X_i}(\frac{X_j}1)=
\begin{cases}
0&\text{if $i\neq j$,}\\
1&\text{otherwise.}
\end{cases}\]

Thus, by adding formal differentiation to $T(\Sigma_m(\overline X))$
we constructed a differential cancellation meadow.

\section{Signed meadows} 
\label{sec:5}
In this section we consider \emph{signed meadows}:
we extend the signature $\Sigma_m=(0,1,+,\cdot,-,^{-1})$ of 
meadows with the unary sign (or signum) function $\sg(x)$.
We write $\Sigma_\textit{ms}$ for this extended signature, so 
$\Sigma_\textit{ms} = (0,1,+,\cdot,-,^{-1},\sg)$.
The sign function $\sg(x)$
presupposes an ordering on its domain and
is defined by

\[\sg(x)=\begin{cases}
-1&\text{if }x<0,\\
0&\text{if }x=0,\\
1&\text{if }x>0.
\end{cases}\]

We define the sign function in an equational manner
by the set \SA\ of axioms 
given in Table~\ref{t:a}.
First, notice that by \Md\ and axiom~\eqref{ax1} 
(or axiom~\eqref{ax2}) we find 
\[\sg(0)=0\quad\text{and}\quad\sg(1)=1.\]
Then, observe that in combination with the inverse law \IL, 
axiom \eqref{ax6} is an equational representation of the conditional
equational axiom
\[\sg(x)=\sg(y)~\longrightarrow ~\sg(x+y)=\sg(x).\]
From \Md\ and axioms
\eqref{ax3}--\eqref{ax6} one 
can easily compute $\sg(t)$ for any closed term $t$.

\begin{table}
\centering
\hrule
\begin{align}
\label{ax1}
\sg(1_x)&=1_x\\
\label{ax2}
\sg(0_x)&=0_x\\
\label{ax3}
\sg(-1)&=-1\\
\label{ax4}
\sg(x^{-1})&=\sg(x)\\
\label{ax5}
\sg(x\cdot y)&=\sg(x)\cdot \sg(y)\\
\label{ax6}
0_{\sg(x)-\sg(y)}\cdot (\sg(x+y)-\sg(x))&=0
\end{align}
\hrule
\caption{The set \SA\ of axioms for the sign function}
\label{t:a}
\end{table}

Some more consequences of the $\Md\cup\SA$ axioms are these:
\begin{align}
\label{1}
\sg(x^2)&=1_x,\\
\label{2}
\sg(x^3)&=\sg(x),\\
\label{3}
1_x\cdot\sg(x)&=\sg(x),\\
\label{5}
\sg(x)^{-1}&=\sg(x).
\end{align}
Here \eqref{1} follows from $\sg(x^2)=\sg(x)\cdot\sg(x)
=\sg(x)\cdot\sg(x^{-1})=\sg(1_x)=1_x$, \eqref{2} from
$\sg(x^3)=\sg(x)\cdot\sg(x)\cdot\sg(x^{-1})
=\sg(x\cdot (x\cdot x^{-1}))=\sg(x)$, \eqref{3} from
$1_x\cdot\sg(x)=\sg(x^2)\cdot\sg(x)
=\sg(x^3)=\sg(x)$, and \eqref{5} from
\begin{align*}
\sg(x)^{-1}&=(\sg(x)^2\cdot\sg(x)^{-1})^{-1}=(\sg(x^2)\cdot\sg(x)^{-1})^{-1}
\\
&=(1_x\cdot\sg(x)^{-1})^{-1}=1_x\cdot \sg(x)=\sg(x).
\end{align*}

So, $0=\sg(x)-\sg(x)=\sg(x)-\sg(x)^3=\sg(x)(1-\sg(x)^2)$ and hence
\begin{equation}
\label{4}
\sg(x)\cdot(1-\sg(x))\cdot (1+\sg(x))=0.
\end{equation}

Identity \eqref{4} implies with \IL\ that for any closed term $t$, 
$\sg(t)\in\{-1,0,1\}$, 
and thus also that $\sg(\sg(t))=\sg(t)$. However, with some 
effort we can derive $\sg(\sg(x))=\sg(x)$, which of course is an
interesting consequence.

\begin{proposition}
\label{Prop}
$\Md\cup\SA\vdash \sg(\sg(x))=\sg(x)$.
\end{proposition}

Before giving a proof of the idempotency of $\sg(x)$
we explain how we found one, as there seems not 
to be an obvious proof for this identity | at the same time
this explanation 
illustrates the proof of Theorem~\ref{st}.
Consider a fresh constant $c$ and let $e$ abbreviate the equation
$\sg(\sg(c)) = \sg(c)$,
then:
\begin{align*}
\Md\cup\SA\cup\{\sg(c) = 0\}
\vdash_\IR e,\\
\Md\cup\SA\cup\{1_{\sg(c)} = 1, 1-\sg(c) = 0\}
\vdash_\IR e,\\
\Md\cup\SA\cup\{1_{\sg(c)} = 1, \;\;1_{1-\sg(c)}=1\} 
\vdash_\IR e.
\end{align*}
The first two derivabilities are trivial, the third one
is obtained from \eqref{4} after multiplication with
$1/\sg(c) \cdot 1/(1-\sg(c))$ (thus yielding $\sg(c)=-1=\sg(\sg(c))$).
The proof transformations that underly the proof of Theorem~\ref{st}
dictate how to eliminate the \IR\ rule in this particular case. 
The proof below shows the slightly polished result.

\begin{proof}[Proof of Proposition~\ref{Prop}.]
Recall $0_t+1_t=1$.
The result $\sg(\sg(x)) = \sg(x)$
follows from 
\begin{align*}
\sg(\sg(x))&=(0_{\sg(x)}+ 1_{\sg(x)})\cdot \sg(\sg(x)),
\\
\sg(x)&=(0_{\sg(x)}+ 1_{\sg(x)})\cdot \sg(x),
\end{align*}
and
\eqref{un} and \eqref{deux}:
\begin{align}
\label{un}
0_{\sg(x)}\cdot \sg(\sg(x))&=0_{\sg(x)}\cdot \sg(x),\\
\label{deux}
1_{\sg(x)}\cdot \sg(\sg(x))&=1_{\sg(x)}\cdot \sg(x).\end{align}
Identity \eqref{un} follows from 
$0=0_{\sg(x)}\cdot \sg(x)$ by 
$0=\sg(0)=\sg(0_{\sg(x)}\cdot \sg(x))
=0_{\sg(x)}\cdot\sg(\sg(x))$, and
\eqref{deux} follows from combining \eqref{trois} and \eqref{quatre}:
\begin{align}
\label{trois}
1_{\sg(x)}\cdot 0_{1-\sg(x)}\cdot\sg(\sg(x))&=
1_{\sg(x)}\cdot 0_{1-\sg(x)}\cdot\sg(x),\\
\label{quatre}
1_{\sg(x)}\cdot 1_{1-\sg(x)}\cdot \sg(\sg(x))&=
1_{\sg(x)}\cdot 1_{1-\sg(x)}\cdot \sg(x).\end{align}
Identity \eqref{trois} follows simply: 
$0_{1-\sg(x)}\cdot(1-\sg(x))=0$,
so $0_{1-\sg(x)}\cdot\sg(x)=0_{1-\sg(x)}$ and thus
\begin{align*}
0_{1-\sg(x)}\cdot\sg(\sg(x))&=\sg(0_{1-\sg(x)}\cdot\sg(x))\\
&=\sg(0_{1-\sg(x)})\\
&=0_{1-\sg(x)}\\
&=0_{1-\sg(x)}\sg(x).
\end{align*}
Identity \eqref{quatre} can be derived as follows: from \eqref{4} infer
\[1_{\sg(x)}\cdot 1_{1-\sg(x)}\cdot (1+\sg(x))=0,\]
thus $1_{\sg(x)}\cdot 1_{1-\sg(x)}\cdot \sg(x)=
1_{\sg(x)}\cdot 1_{1-\sg(x)}\cdot -1$, and thus with $\sg(-1)=-1$,
\begin{align*}
1_{\sg(x)}\cdot 1_{1-\sg(x)}\cdot\sg(\sg(x))&= 
\sg(1_{\sg(x)}\cdot 1_{1-\sg(x)}\cdot \sg(x))\\
&=
1_{\sg(x)}\cdot 1_{1-\sg(x)}\cdot -1\\
&=1_{\sg(x)}\cdot 1_{1-\sg(x)}\cdot \sg(x).
\end{align*}
\end{proof}

Next we establish the expected corollary of Theorem~\ref{st}:
\begin{corollary}
\label{cor:signed}
The set of axioms $\Md\cup\SA$ (see Tables~\ref{Md} and \ref{t:a})
is a finite basis (a complete axiomatization) of 
$\Mod_{\Sigma_{ms}}(\Md\cup\SA\cup\IL)$.
\end{corollary}
\begin{proof}
It suffices to show that
the propagation properties are satisfied for $\sg(\_)$.
\bigskip\\
Pseudo units:
$1_x\cdot\sg(y)=(1_x)^2\cdot\sg(y) =1_x\cdot\sg(1_x)\cdot \sg(y)
=1_x\cdot\sg(1_x\cdot y)$.
\bigskip\\
Pseudo zeros:
$0_x\cdot\sg(y)=(0_x)^2\cdot\sg(y) =0_x\cdot\sg(0_x)\cdot\sg(y)
=0_x\cdot\sg(0_x\cdot y)$.
\end{proof}

We notice that the initial algebra
of $\Md \cup \SA$ equals $\mathbb{Q}_0$ as introduced in~\cite{BT07}
expanded with the sign function
(a proof follows immediately from the techniques used in that paper).
It remains to be shown that the \SA\ axioms (in combination with 
those of \Md)
are independent. We leave this as an open question.

In the following we show that the sign function is not
definable in $\mathbb{Q}_0$, the zero-totalized field 
of rational numbers as discussed in~\cite{BT07}.
 We say that
$q,q'\in T(\mathbb{Q}_0)$ are \emph{different}
if $1_{q-q'}=1$. Let $r=r(x)$ and $s=s(x)$ and let $T(\mathbb{Q}_0[x])$
be the set of terms that are either closed or have $x$ as
the only variable,
so $r,s\in T(\mathbb{Q}_0[x])$. We define
\begin{eqnarray*}
r\equiv_\infty s&\iff&r(q)=s(q) \quad
\begin{array}[t]{l}
\text{for infinitely many}\\
\text{different $q$ in $T(\mathbb{Q}_0)$},\end{array}\\
r\equiv_{ae} s&\iff&r(q)\neq s(q) \quad
\begin{array}[t]{l}
\text{for finitely many}\\
\text{different $q$ in $T(\mathbb{Q}_0)$}.\end{array}
\end{eqnarray*}
We call these relations \emph{infinite 
equivalence} and \emph{almost equivalence}, respectively.
Observe that both these relations are congruences over $T(\mathbb{Q}_0[x])$. 

\begin{theorem}
\label{ae}
Let $r=r(x)$ and $s=s(x)$. If $r\equiv_\infty s$ then 
$r\equiv_{ae} s$.
\end{theorem}

\begin{proof}
By Theorem~\ref{SMF} it
suffices to prove this for SMFs, say $P=P(x)$ and $Q=Q(x)$.  
Because $P-Q$ is then provably equal to an SMF, we 
further assume without loss of generality that  $Q= 0$. 

So, let $P\equiv_\infty 0$.
We prove $P\equiv_{ae} 0$ by induction on the level $n$ of $P$.
\begin{description}
\item[Case $n=0$.] 
Then $P=s/t$ for polynomials $s=s(x)$ and $t=t(x)$. 
Because $P\equiv_\infty 0$, at least one of 
$s\equiv_\infty 0$ and $t\equiv_\infty 0$ holds.
Because polynomials always have a finite number of zero points, 
at least one of 
$s\equiv_{ae} 0$ and $t\equiv_{ae} 0$ holds. Thus $P\equiv_{ae} 0$.
\item[Case $n+1$.] 
Then $P=0_t\cdot S+1_t\cdot T$.
\begin{itemize}
\item If $t\equiv_{ae}0$ then $0_t\equiv_{ae}1$ and 
$1_t\cdot T\equiv_{ae}0$, so $S\equiv_\infty 0$. By induction,
$S\equiv_{ae}0$, and thus $0_t\cdot S\equiv_{ae}0$ and hence 
$P\equiv_{ae} 0$.
\item If $t\not\equiv_{ae}0$ then $1_t\equiv_\infty 1$, so 
$1_t\equiv_{ae}1$ and $0_t\cdot S\equiv_{ae}0$, so $T\equiv_\infty 0$.
By induction,
$T\equiv_{ae}0$, and thus $1_t\cdot T\equiv_{ae}0$ and hence 
$P\equiv_{ae} 0$.

\end{itemize}

\end{description}
\end{proof}

An immediate consequence of Theorem~\ref{ae} is:
\begin{corollary}
The sign function is not definable in $\mathbb{Q}_0$.
\end{corollary}
\begin{proof} Suppose otherwise. 
Then there is a term $t\in T(\mathbb{Q}_0[x])$ with 
$\sg(x)=t(x)$. So 
\[t(x)\equiv_\infty 1\]
(because of all positive rationals). But then $t(x)
\equiv_{ae}1$ by Theorem~\ref{ae}, which contradicts $t(x)=-1$ 
for all negative rationals.
\end{proof}

Furthermore, we notice that with the sign function $\sg(x)$,
the functions $\max(x,y)$ and $\min(x,y)$ have a simple equational specification:
\begin{align*}
\max(x,y)&=\max(x-y,0)+y,\\
\max(x,0)&=(\sg(x)+1)\cdot x/2,
\end{align*}
and, of course, $\min(x,y)=-\max(-x,-y)$.

Finally, the existence of
non-trivial differential cancellation meadows with sign function is not an 
obvious matter and requires a modification of the 
existence proof given in~Section~\ref{subsec:4.2}.

\section{Floor, Ceiling and Square Root}
\label{sec:6}
In this section we consider extensions of signed meadows with floor, ceiling and square root.

\subsection{Signed Meadows with Floor and Ceiling}
We briefly discuss the extension of
signed meadows with the \emph{floor} function $\floor x$
and the \emph{ceiling} function $\ceiling x$. These functions
are defined by
\[\floor x =\max\{n\in\NZ\mid n\leq x\}\]
and 
\[\ceiling x=\min\{n\in\NZ\mid n\geq x\}.
\]
We define these functions in an equational manner by the axioms
in Table~\ref{floor}.

\begin{table}
\centering
\hrule
\begin{align}
\label{floor1}
1_x\cdot \floor{y}&=1_x\cdot\floor{1_x\cdot y}\\
\label{floor2}
0_x \cdot \lfloor y \rfloor &=  
0_x \cdot \lfloor 0_x \cdot y \rfloor\\
\label{floor3}
\floor {x - 1}&=\floor x -1
\\
\label{floor4}
\floor {x + 1}&=\floor x +1
\\
\label{floor5}
\floor {0}&=0
\\
\label{floor6}
(0_{1-\sg(x)}\cdot0_{1-\sg(1-x)})\cdot\floor x&=0
\\[2mm]
\label{ceil1}
\ceiling x&= -\floor{-x}
\end{align}
\hrule
\caption{The set \FC\ of axioms for the floor and ceiling functions}
\label{floor}
\end{table}

Some comments on these axioms: first, 
\eqref{floor1} and \eqref{floor2}
guarantee the propagation properties. Then, consider
$0_{1-\sg(x)}\cdot 0_{1-\sg(1-x)}$, which equals 1
if both $x>0$ and $1-x>0$, and 0 otherwise. 
So, axiom \eqref{floor6} states 
that $\floor{x}=0$ whenever $0<x<1$. 
With \eqref{floor3}--\eqref{floor5} this is sufficient to
compute $\floor t$ for any closed $t$. Axiom~\eqref{ceil1}, 
defining the ceiling function $\ceiling x$ is 
totally standard.

Let $\Sigma_\textit{msfc}$ be the signature of this extension.
As before, we have an immediate corollary of Theorem~\ref{st}.
\begin{corollary}
The set of axioms $\Md\cup\SA\cup\FC$ (see Tables~\ref{Md},
\ref{t:a} and \ref{floor})
is a finite basis (a complete axiomatization) of 
$\Mod_{\Sigma_\textit{msfc}}(\Md\cup\SA\cup\FC\cup\IL)$.
\end{corollary}
\begin{proof}
For floor, the propagation properties for pseudo units and 
for pseudo zeros are directly axiomatized by axioms \eqref{floor1} 
and \eqref{floor2}, and those for ceiling follow easily. 
So, the 
corollary follows immediately from Theorem~\ref{st} and the 
proof of Corollary~\ref{cor:signed}.
\end{proof}

We notice that the initial algebra of $\Md\cup\SA\cup\FC$ is $\mathbb{Q}_0$ 
extended with the sign function $\sg(x)$ and the floor and ceiling
functions $\floor x$ and $\ceiling x$.
It remains to be shown that the \FC\ axioms (in combination with 
those of $\Md\cup\SA$)
are independent. We leave this as an open question.

We continue this section by proving that in $\mathbb{Q}_0(\sg)$, i.e.,
the rational numbers viewed as a signed meadow, a 
definition of ceiling and floor cannot be given. To this end,
we first prove a general property of unary functions definable in 
$\mathbb{Q}_0(\sg)$.

\begin{theorem}
\label{explicit}
For any function $h(x)$ definable in 
$\mathbb{Q}_0(\sg)$ there exist $r\in T(\mathbb{Q}_0)$ and a function $g(x)$ 
definable in $\mathbb{Q}_0[x]$ such that
\[
x>r~\Longrightarrow~ h(x)=g(x).\]
\end{theorem}

\begin{proof}
By structural induction on the form that $h(x)$ may take. 

If $h(x)\in\{0,1, x\}$, we're done. For $h(x)=-f(x)$ or 
$h(x)=1/f(x)$ or
$h(x)=f_1(x)+f_2(x)$ or $h(x)=f_1(x)\cdot f_2(x)$, the result 
also follows immediately
(in the latter cases take $r=\max(r_1,r_2)$ for $r_i$
satisfying the property for $f_i(x)$).

In the remaining case,
$h(x)=\sg(f(x))$. Let 
$g(x)\in T(\mathbb{Q}_0[x])$ 
be such that $f(x)=g(x)$ for $x>r$.
By induction on the form that $g(x)$ may take, it follows
that an $r'$ exists
such that for $x>r'$, $\sg(g(x))$ is constant. This proves 
that for $x>\max(r,r')$, $h(x)=\sg(f(x))=\sg(g(x))$ is constant.
\end{proof}

\begin{corollary} The floor function $\floor x$ is not 
definable in $\mathbb{Q}_0(\sg)$.
\end{corollary}

\begin{proof} Consider
\[h(x) = \frac{x-\floor x}{x-\floor x}.\]
If $h(x)$ were definable in $\mathbb{Q}_0(\sg)$, 
then by the preceding result 
there exist $r$ and a function $g(x)$ definable in $\mathbb{Q}_0[x]$
such that $h(x)=g(x)$ for $x>r$. But then 
$g(x)\equiv_{\infty}0$ (for all integers above $r$)
and $g(x)\equiv_\infty 1$ (for all non-integers above $r$),
and this contradicts Theorem~\ref{ae}.
\end{proof}

We finally notice that for $t(x)$ some term one can add this 
induction rule:
\[
\frac{\begin{array}{l}t(0)=0, \\
0_{1-\sg(x)}\cdot0_{t(\floor x)}\cdot t(\floor x + 1)=0,\\
0_{1+\sg(x)}\cdot0_{t(\ceiling x)}\cdot t(\ceiling x - 1)=0
\end{array}}
{t(\floor x)=0,\quad t(\ceiling x)=0}
\]
thus
\[
\frac{\begin{array}{l}t(0)=0, \\
(x>0~\&~t(\floor x)=0) ~\longrightarrow~ t(\floor x + 1)=0,\\
(x<0~\&~t(\ceiling x)=0)~\longrightarrow~ t(\ceiling x - 1)=0
\end{array}}
{t(\floor x)=0,\quad t(\ceiling x)=0}.
\]
With this particular induction rule,
the idempotency of $\floor x$ can be easily proved 
(take $t(x)=x-\floor x$), as well as the idempotency of
ceiling. With a little more effort one can prove 
$\floor{ x -\floor x} = 0$:
first prove $ \floor{-\floor x}=-\floor x$ by induction on $x$,
and then
$\floor{ x + \floor y} = \floor x + \floor y$ by induction on $y$. 
As a consequence, 
$\floor{ x - \floor x}= \floor x + \floor{-\floor x} 
= \floor x + -\floor x= 0$.
In general, if using \IL\ the premises
can be proved (from some extension of \Md\ that satisfies 
the propagation properties), then this can also be proved 
without \IL,
and therefore this also is the case for the conclusion.

\subsection{Signed Meadows with Square Root}
A plausible way to totalize the square root operation is to 
postulate $\sqrt{-1}=i$ and to abandon the 
domain of signed fields in favour of the complex numbers. 
Here we choose a different approach by stipulating
$\sqrt{x}=-\sqrt{-x}$ for $x<0$. In order to avoid confusion with 
the principal square root function we deviate from the standard
notation and introduce the unary operation $\wortel{\_}$ called 
\emph{signed square root}.
We write $\Sigma_{mss}$ for this extended signature, so 
$\Sigma_{mss} = (0,1,+,\cdot,-,^{-1},\sg,\  \wortel{\ })$,
and define the signed square root operation in an equational manner
by the set \SR\ of axioms 
given in Table~\ref{t:wortel}.

\begin{table}
\centering
\hrule
\begin{align}
\label{ax7}
\wortel{x^{-1}}&=(\wortel{x})^{-1}\\
\label{ax8}
\wortel{x\cdot y}&=\wortel{x}\cdot \wortel{y}\\
\label{ax9}
\wortel{x\cdot x\cdot \sg(x)}&=x\\
\label{ax10}
\sg(\wortel{x}-\wortel{y})&=\sg(x-y)
\end{align}
\hrule
\caption{The set \SR\ of axioms for the square root}
\label{t:wortel}
\end{table}
\noindent

Some additional consequences of the $\Md\cup\SA\cup \SR$ axioms are these:
\begin{align}
\label{6}
\wortel{\sg(x)}&=\sg(x),\\
\label{7}
\wortel{1_x}&=1_x,\\
\label{7a}
\wortel{0_x}&=0_x,\\
\label{8}
\wortel{-x}&=-\wortel{x},\\
\label{9}
\wortel{x^2}&=x\cdot \sg(x).
\end{align}
Here identity \eqref{6} follows from 
\begin{align*}
\wortel{\sg(x)}&=\wortel{\sg(xxx^{-1})}\\
&=\wortel{\sg(x)\sg(x)\sg(x^{-1})}\\
&=\wortel{\sg(x)\sg(x)\sg(x)}\\
&=\wortel{\sg(x)\sg(x)\sg(\sg(x))}=\sg(x),
\end{align*}
identity~\eqref{7} from $\wortel{1_x}=\wortel{\sg(1_x)}=\sg(1_x)=1_x$ 
and identity \eqref{7a} is proved similarly.
Identity~\eqref{8} follows from
\begin{align*}
\wortel{-x}&=\wortel{-1\cdot x}
=\wortel{-1}\cdot \wortel{x}
=\wortel{\sg(-1)}\cdot\wortel{x}\\
&=\sg(-1)\cdot\wortel{x}=-1\cdot \wortel x=-\wortel x,
\end{align*} 
and \eqref{9} from
\begin{align*}
\wortel{x^2}&=\wortel{x^2\cdot 1_x}
=\wortel{x^2}\cdot 1_x=\wortel{x^2}\cdot \sg(1_x)\\
&=\wortel{x^2}\cdot \sg(x)^2
=\wortel{x^2}\cdot \wortel{\sg(x)}\cdot \sg(x)\\
&=\wortel{x^2\sg(x)}\cdot \sg(x)=
x\cdot \sg(x).
\end{align*}

Since $(\Sigma_{mss},\Md \cup \SA \cup \SR)$ satisfies both propagation properties, 
we can apply Theorem~\ref{st}. 
\begin{corollary}
The set of axioms $\Md\cup \SA\cup \SR$
is a complete axiomatization of 
$\Mod_{\Sigma_{mss}}(\Md\cup \SA \cup \SR \cup\IL)$. 
\end{corollary}
\begin{proof}
We have to prove that the propagation properties for pseudo units 
and pseudo zeros hold in $\Md\cup \SA \cup \SR$. This follows easily
by a case distinction on the forms that $C[r]$
may take. 
As an example we consider here the case
$C[\_]\equiv \wortel{\_}$. Then 
\[
1_t\cdot \wortel{r}=1_t^2\cdot \wortel{r}=1_t\cdot \wortel{1_t} \cdot \wortel{r}=
1_t\cdot \wortel{1_t\cdot r}
\]
by~\eqref{eq:1} and~\eqref{7}.
The propagation property for pseudo zeros is proved in a similar way 
applying~\eqref{eq:2} and~\eqref{7a}.
\end{proof}

We denote by $\mathbb{Q}_0(\sg , \wortel{\ })$ the zero-totalized signed prime field that contains
$\mathbb{Q}$ and is closed under $\wortel{\ }$. 
Note that $\mathbb{Q}_0(\sg , \wortel{\ })$ 
is a computable data type (see e.g. Bergstra and 
Tucker~\cite{BT95}). This statement still requires an efficient 
and readable proof.

Finally, differential meadows can be equipped with a signed 
square root operator by the axioms given in Table~\ref{t:diff}. 
Axiom~\eqref{ax15} can actually be derived from Axiom~\eqref{ax14} and the 
equational axiomatization of differential meadows as follows:

\begin{align*}
2\cdot \wortel{y}\cdot\frac{\partial}{\partial X_i}(\wortel{y})
&=\wortel{y}\cdot\frac{\partial}{\partial X_i}(\wortel{y}) + 
\wortel{y}\cdot\frac{\partial}{\partial X_i}(\wortel{y})\\
&\stackrel{\eqref{D2}}=\frac{\partial}{\partial X_i}(\wortel{y}\cdot \wortel{y})\\
&\stackrel{\eqref{ax8}}=\frac{\partial}{\partial X_i}(\wortel{y^2})\\
&\stackrel{\eqref{9}}=\frac{\partial}{\partial X_i}(y\cdot \sg(y))\\
&\stackrel{\eqref{D2}}=\sg(y)\cdot \frac{\partial}{\partial X_i}(y) + 
y\cdot \frac{\partial}{\partial X_i}(\sg(y)) \\
&\stackrel{\eqref{ax14}}=\sg(y)\cdot \frac{\partial}{\partial X_i}(y).\\
\end{align*}
Moreover,  by identity~\eqref{7}, $1_y = 1_{\wortel{y}}$, and thus 
\[\wortel{y} = \wortel{1_y \cdot y} = 1_y \cdot \wortel{y}.\]
Hence
\begin{align*}
\frac{\partial}{\partial X_i}(\wortel{y})
&\stackrel{\phantom{\eqref{7}}}= \frac{\partial}{\partial X_i}(1_y \cdot \wortel{y})\\
&\stackrel{\eqref{D2}}=\wortel{y}\cdot\frac{\partial}{\partial X_i}(1_y )+
1_y\cdot\frac{\partial}{\partial X_i}(\wortel{y})\\
&\stackrel{\eqref{D3}}=1_y\cdot\frac{\partial}{\partial X_i}(\wortel{y})\\
&\stackrel{\eqref{7}}=1_{\wortel{y}} \cdot  \frac{\partial}{\partial X_i}(\wortel{y})\\
&\stackrel{\phantom{\eqref{7}}}=\frac{\sg(y)}{2}(\wortel{y})^{-1}\cdot \frac{\partial}{\partial X_i}y.
\end{align*}
So, the existence of
non-trivial differential cancellation meadows with signed square roots depends heavily 
on the existence of an appropriate interpretation of the sign function.

\begin{table}
\centering
\hrule
\begin{align}
\label{ax14}
\frac{\partial}{\partial X_i}\sg(y)&=0\\
\label{ax15}
\frac{\partial}{\partial X_i}\wortel{y}
&=\frac{\sg(y)}{2}(\wortel{y})^{-1}\cdot \frac{\partial}{\partial X_i}y
\end{align}
\hrule
\caption{The signed square root for differential meadows}
\label{t:diff}
\end{table}
\noindent

\section{Conclusions}
\label{sec:7}
The main result of this paper is a generic  basis theorem
for cancellation meadows. We have applied this result 
to various expansions of meadows. The first expansion concerns differential fields.  
It appears that the 
interaction between differential operators and equations for 
meadows is entirely unproblematic. The propagation properties 
follow immediately from well-known axioms for differential 
fields.

As stated before, most uses of rational numbers 
in computer science exploit their ordering. We 
include this ordering by extending the initial algebraic 
specification of $\mathbb{Q}_0$ with an equational specification of the 
sign function, resulting in a finite basis for what we
called $\mathbb{Q}_0(\sg)$ and we provided a non-trivial proof of the 
idempotency of the sign function in $\mathbb{Q}_0(\sg)$. 
However, the question whether our
particular axioms for $\sg(x)$ are independent is left open.

As a further example we added the floor function $\floor x$, 
the ceiling function $\ceiling x$, and the signed square root  to 
$\mathbb{Q}_0(\sg)$ and showed that the resulting equational specification 
is a finite basis. 
Again, we did not investigate the independency
of these axioms.

In~\cite{BT95} it is shown that computable algebras can be 
specified 
by means of a complete term rewrite system, provided auxiliary 
functions can be used. Useful candidates for auxiliary operators in
the case of rational numbers can be found in Moss~\cite{MOSS} and 
Calkin and Wilf~\cite{CW}.
In~\cite{BT07} the existence of an equational 
specification of $\mathbb{Q}_0$ which is confluent and terminating as a rewrite 
system has been formulated as an open question. 
To that question we now add the 
corresponding question in the presence of the sign operator.

\end{document}